\documentclass[12pt,a4paper,fleqn]{article}
\usepackage{a4wide,amsfonts,amsmath,latexsym,amssymb,euscript,graphicx,units,mathrsfs}
\usepackage{mathpazo}
\usepackage{hyperref}
\usepackage{multimedia}

\newtheorem{theorem}{Theorem}

\newtheorem{example}[theorem]{Example}

\newtheorem{proposition}[theorem]{Proposition}
\newtheorem{remark}[theorem]{Remark}

\newenvironment{proof}[1][Proof]{\textbf{#1.} }{\ \rule{0.5em}{0.5em}}

\begin{document}

\title{Asymptotics of weighted \ random sums }
\author{Jos\'{e} Manuel Corcuera\thanks{%
University of Barcelona, Spain. \texttt{E-mail: jmcorcuera@ub.edu}}, David
Nualart\thanks{%
University of Kansas, USA. \texttt{E-mail: nualart@ku.edu}}, Mark Podolskij%
\thanks{%
University of Heilderberg, Germany. \texttt{E-mail:}%
m.podolskij@uni-heidelberg.de}}
\maketitle

\begin{abstract}
In this paper we study the asymptotic behaviour of weighted random sums when
the sum process converges stably in law to a Brownian motion and the weight
process has continuous trajectories, more regular than that of a Brownian
motion. We show that these sums converge in law to the integral of the
weight process with respect to the Brownian motion when the distance between
observations goes to zero. The result is obtained with the help of
fractional calculus showing the power of this technique. This study, though
interesting by itself, is motivated by an error found in the proof of
Theorem 4 in \cite{CorNuaWoe06}.
\end{abstract}

\section{Introduction}

Let $\left( \Omega ,\mathcal{F}, \mathbb{P}\right) $ be a complete
probability space. Fix a time interval $[0,T]$ and consider a double
sequence of random variables $\xi=\{\xi_{i,m}, m \in \mathbb{Z}_+, 1\le i
\le [mT]\}$. For any $m\ge 1$ we denote by $g_m(t)$ the stochastic process
defined as the distribution function of the signed measure on $[0,T]$ which
gives mass $\xi_{i,m}$ to the points $t_i= \frac im$, for $1\le i \le [mT]$,
that is,
\begin{equation*}
g_m(t):=\sum_{i=1}^{[mt]} \xi_{i,m}.
\end{equation*}%
Notice that $\xi_{i,m}= g_m(t_i)- g_m(t_{i-1})$. Throughout this paper we
assume the following hypotheses on the double sequence $\xi$:

\medskip \noindent\textbf{(H1)} The sequence of processes $\left(
g_m(t)\right) _{t\in [0,T]}$ satisfies
\begin{equation*}
\left( g_m(t)\right) _{t\in [0,T]}\overset{f.d.d.}{\underset{m\rightarrow
\infty }{\rightarrow }}\left( w(t)\right)_{t\in [0,T]} \qquad \mathcal{F}%
\text{-stably},
\end{equation*}%
where $\left( w(t)\right) _{t\in [0,T]}$ is a standard Brownian motion
independent of $\mathcal{F}$, and the latter denotes convergence of finite
dimensional distributions $\mathcal{F}$-stably in law (see the definition
below).

\medskip \noindent\textbf{(H2)} The family of random variables $\xi$
satisfies the tightness condition
\begin{equation}
\mathbb{E} \left( \left\vert \sum_{i=j+1}^{k} \xi_{i,m} \right\vert
^{4}\right) \leq C\left( \frac{k-j}{m}\right) ^{2},  \label{tight}
\end{equation}%
for any $1\le j <k \le [mT]$.

\medskip Notice that these hypotheses imply that $\left( g_{m}(t)\right)
_{t\in \lbrack 0,T]}\overset{}{\overset{\mathcal{L}}{\underset{m\rightarrow
\infty }{\rightarrow }}}\left( w(t)\right) _{t\in \lbrack 0,T]}$, $\mathcal{F%
}\text{-stably}$ in the Skorohod space $D[0,T]$ equipped with the uniform
topology.

\medskip Under these assumptions, the purpose of this note is to establish
the following result.

\begin{theorem}
\label{thm1} Let $\left( f(t)\right) _{t\in [0,T]}$ be a $\alpha $-H\"{o}%
lder continuous process with index $\alpha >1/2.$ Suppose that $%
\xi=\{\xi_{i,m}, m \in \mathbb{Z}_+, 1\le i \le [mT]\}$ is a family of
random variables satisfying hypotheses \textbf{(H1)} and \textbf{(H2)}. Set
\begin{equation*}
X_m(t):=\sum_{i=1}^{[mt]}f(t_{i}) \xi_{i,m}=\int_{0}^{t}f(s)\mathrm{d}%
g_{m}(s).
\end{equation*}
Then,
\begin{equation*}
X_m(t)\overset{\mathcal{L}}{\underset{m\rightarrow \infty }{\rightarrow }}%
\int_{0}^{t}f(s)\mathrm{d}w(s), \qquad \mathcal{F}\text{-stably}
\end{equation*}
in $D[0,T]$, where $\left( w(t)\right) _{t\in [0,T]}$ is a standard Brownian
motion independent of $\mathcal{F}$.
\end{theorem}

Recall that a sequence of random vectors or processes $Y_n$ converges $%
\mathcal{F}$-stably in law to a random vector or process $Y$, where $Y$ is
defined on an extension $(\Omega^{\prime }, \mathcal{F}^{\prime }, \mathbb{P}%
^{\prime })$ of the original probability $(\Omega, \mathcal{F}, \mathbb{P})$%
, if $\ (Y_{n},Z)\overset{\mathcal{L}}{\rightarrow }(Y,Z)$ for any $\mathcal{%
F}$-measurable random variable $Z$. If $Y$ is $\mathcal{F}$-measurable, then
we have convergence in probability. We refer to \cite{Ren63} and \cite%
{AldEag78} for more details on stable convergence.

\begin{remark}
The conclusion of Theorem \ref{thm1} also holds for the forward-type Riemann
sums
\begin{equation*}
\widetilde{X}_m(t):=\sum_{i=1}^{[mt]}f(t_{i-1}) \xi_{i,m}.
\end{equation*}
Indeed, it suffices to put the random weights $\xi_{i,m}$ at the points $%
t_{i-1}$.
\end{remark}

This theorem has been motivated by a mistake in the proof of Theorem 4 of
the reference \cite{CorNuaWoe06}. More precisely the fact that $%
\lim_{n\rightarrow \infty }\limsup_{m\rightarrow \infty }\mathbb{P}(\Vert
B_{t}^{(n,m)}\Vert _{\infty }>\epsilon )=0$ in page 724 of \cite{CorNuaWoe06}
is a particular case of the convergence (\ref{vancov}). The rest of this
note is devoted to the proof of Theorem \ref{thm1}. First we establish a
basic decomposition, which reduces the proof of Theorem \ref{thm1} to the
proof of the convergence (\ref{vancov}). In section 3 we discuss different
attempts to prove this convergence using $p$-variation norms and martingale
methods. Finally in Section 4 we provide a proof of (\ref{vancov}) using
techniques of fractional calculus.

\section{The main decomposition}

The basic idea of the proof of Theorem \ref{thm1} is the classical
Bernstein's big blocks/small blocks technique. For this purpose we set $%
u_{j}=\frac{j}{n},n\leq m$, and decompose the process $X_{m}(t)$ as follows
\begin{eqnarray}
X_{m}(t) &=&\sum_{i=1}^{[mt]}f(t_{i})\xi _{i,m}  \notag \\
&=&\sum_{j=1}^{[nt]+1}\sum_{i\in I_{n}(j)}^{{}}\left(
f(t_{i})-f(u_{j-1})\right) \xi
_{i,m}+\sum_{j=1}^{[nt]+1}f(u_{j-1})\sum_{i\in I_{n}(j)}^{{}}\xi _{i,m}
\label{block}
\end{eqnarray}%
with $I_{n}(j):=\{i:1\leq i\leq \lbrack mT],\frac{i}{m}\in \lbrack \frac{j-1%
}{n},\frac{j}{n})\}.$ \ According to our hypothesis it holds that
\begin{equation*}
g_{m}\overset{}{\overset{\mathcal{L}}{\underset{m\rightarrow \infty }{%
\rightarrow }}}w\qquad \mathcal{F}\text{-stably}
\end{equation*}%
on $D[0,T]$ with $\left( w(t)\right) _{t\in \lbrack 0,T]}$ being a Brownian
motion independent of $\mathcal{F}$. This implies, in particular, the $%
\mathcal{F}$-stable convergence
\begin{equation*}
\sum_{j=1}^{[nt]+1}f(u_{j-1})\sum_{i\in I_{n}(j)}^{{}}\xi _{i,m}\overset{}{%
\overset{\mathcal{L}}{\underset{m\rightarrow \infty }{\rightarrow }}}%
\sum_{j=1}^{[nt]+1}f(u_{j-1})\left( w(u_{j})-w(u_{j-1})\right) .
\end{equation*}%
We also have that
\begin{equation*}
\sum_{j=1}^{[nt]+1}f(u_{j-1})\left( w(u_{j})-w(u_{j-1})\right) \overset{}{%
\overset{u.c.p.}{\underset{n\rightarrow \infty }{\rightarrow }}}%
\int_{0}^{t}f(s)\mathrm{d}w(s).
\end{equation*}%
where $u.c.p.$ stands for uniform convergence in probability.

Now we treat the first term of \ (\ref{block}), but before we consider,
separately, the last summand. We claim that
\begin{equation}
\mathbb{P}\text{-}\lim_{n\rightarrow \infty }\limsup_{m\rightarrow \infty
}\sup_{t\in \lbrack 0,T]}\left\vert \sum_{i\in I_{n}([nt]+1)}\left(
f(t_{i})-f(u_{[nt]})\right) \xi _{i,m}\right\vert =0.  \label{eq1}
\end{equation}%
In fact, using the H\"{o}lder continuity of $f$ we can write
\begin{equation*}
\left\vert \sum_{i\in I_{n}([nt]+1)}\left( f(t_{i})-f(u_{[nt]})\right) \xi
_{i,m}\right\vert \leq \Vert f\Vert _{\alpha }n^{-\alpha }\sum_{i\in
I_{n}([nt]+1)}|\xi _{i,m}|,
\end{equation*}%
where
\begin{equation*}
||f||_{\alpha }:=\sup_{|u-v|\leq T}\frac{|f(u)-f(v)|}{|u-v|^{\alpha }}%
<\infty .
\end{equation*}%
Then, (\ref{eq1}) follows from Hypothesis \textbf{(H2)} taking into account
that the cardinality of the set $I_{n}([nt]+1)$ is bounded by $\frac{m}{n}+1$
and $\alpha >\frac{1}{2}$.

Then, in order to finish the proof, we need to show that
\begin{equation}
\mathbb{P}\text{-}\lim_{n\rightarrow \infty }\limsup_{m\rightarrow \infty }
\sup_{t\in \lbrack 0,T]} \left| \sum_{j=1}^{[nt]}\sum_{i\in
I_{n}(j)}^{{}}\left( f(t_{i})-f(u_{j-1})\right) \xi_{i,m} \right|=0.
\label{vancov}
\end{equation}

In fact, this is a key step of the proof. In particular situations, such as
e.g. in the martingale framework, there are various specific techniques of
the proof. We will present some of them in the next section. However,
proving convergence \eqref{vancov} in a general setting turns out to be not
quite easy.

A first straightforward attempt is as follows. We set
\begin{equation*}
R_{n,m}(t):=\sum_{j=1}^{[nt]}\sum_{i\in I_{n}(j)}^{{}}\left(
f(t_{i})-f(u_{j-1})\right) \xi _{i,m}.
\end{equation*}%
Then we deduce
\begin{equation*}
\sup_{t\in \lbrack 0,T]}\left\vert R_{n,m}(t)\right\vert \leq n^{-\alpha
}||f||_{\alpha }\sum_{i=1}^{[mT]}\left\vert \xi _{i,m}\right\vert ,
\end{equation*}%
but in general we have that
\begin{equation*}
\lim_{m\rightarrow \infty }\sum_{i=1}^{[mT]}\left\vert \xi _{i,m}\right\vert
=\infty
\end{equation*}%
as the following simple example shows.

\begin{example}
\textrm{Consider the case where $\xi_{i,m}= \frac {X_i}{ \sqrt{m}}$, where $%
\left( X_{i}\right) _{i\geq 1}$ are i.i.d with $\mathbb{P}(X _{i}=1)=\mathbb{%
P}(X_{i}=-1)=1/2.$ Then
\begin{equation*}
g_m(t) =\sum_{i=1}^{[mt]}\frac{X_{i}}{\sqrt{m}},
\end{equation*}%
and $g_m\overset{\mathcal{L}}{\underset{m\rightarrow \infty }{\rightarrow }}%
w $ on $D[0,T]$ and
\begin{equation*}
\sum_{i=1}^{[mT]}\left\vert \xi_{i,m}\right\vert \underset{m\rightarrow
\infty }{\rightarrow }\infty.
\end{equation*}
}
\end{example}

\section{Young's calculus}

A more sophisticated approach for proving \eqref{vancov} is to use Young's
integral. Consider the interval $I_{j}^{n}:=\left[ u_{j-1},u_{j}\right) $.
Then
\begin{equation*}
\sum_{j=1}^{[nt]}\sum_{i\in I_{n}(j)}\left( f(t_{i})-f(u_{j-1})\right) \xi
_{i,m}=\sum_{j=1}^{[nt]}\int_{I_{j}^{n}}(f(s)-f(u_{j-1}))\mathrm{d}g_{m}(s).
\end{equation*}%
By the Love-Young inequality and for $\beta >1-\alpha ,$
\begin{equation*}
\left\vert \sum_{j=1}^{[nt]}\int_{I_{j}^{n}}(f(s)-f(u_{j-1}))\mathrm{d}%
g_{m}(s)\right\vert \leq C_{\alpha ,\beta }\sum_{j=1}^{[nt]}\upsilon _{\frac{%
1}{\alpha }}^{j}(f)\upsilon _{\frac{1}{\beta }}^{j}(g_{m})
\end{equation*}%
with
\begin{equation*}
\upsilon _{p}^{j}(h):=\left( \sup_{\pi }\sum_{i=1}^{N}\left\vert
h(s_{i})-h(s_{i-1})\right\vert ^{p}\right) ^{1/p},
\end{equation*}%
where the supremum runs over all partitions $\pi =\{s_{0},\dots ,s_{N}\}$ of
the interval $[u_{j-1},u_{j}]$, and $C_{\alpha ,\beta }$ is a positive
constant; see \cite{You36}. Notice that $\upsilon _{\frac{1}{\alpha }%
}^{j}(f)\leq n^{-\alpha }||f||_{\alpha }$. The problem is then to bound $%
\upsilon _{\frac{1}{\beta }}^{j}(g_{m})$ under the hypothesis $g_{m}\overset{%
\mathcal{L}}{\underset{m\rightarrow \infty }{\rightarrow }}w$ on $D[0,T]$.
Unfortunately, the strong $p$-variation $\upsilon _{p}$ is not a continuous
functional on $D[0,T]$ equipped with the uniform topology. Nevertheless, we
suppose for a moment that the convergence
\begin{equation*}
\upsilon _{\frac{1}{\beta }}^{j}(g_{m})\overset{\mathcal{L}}{\underset{%
m\rightarrow \infty }{\rightarrow }}\upsilon _{\frac{1}{\beta }}^{j}(w),
\end{equation*}%
holds. Then, recalling that $\alpha >\frac{1}{2}$, we can choose $\beta =%
\frac{1}{2}-\varepsilon $ with $\varepsilon <\alpha -\frac{1}{2}$ and we
obtain that
\begin{equation*}
\upsilon _{\frac{1}{\beta }}^{j}(w)\leq ||w||_{\beta }n^{-\beta }<\infty .
\end{equation*}%
Consequently,
\begin{equation*}
\upsilon _{\frac{1}{\alpha }}^{j}(f)\upsilon _{\frac{1}{\beta }}^{j}(g_{m})%
\overset{\mathcal{L}}{\underset{m\rightarrow \infty }{\rightarrow }}\upsilon
_{\frac{1}{\alpha }}^{j}(f)\upsilon _{\frac{1}{\beta }}^{j}(w)\leq
||f||_{\alpha }||w||_{\beta }n^{-\alpha -\beta }.
\end{equation*}%
Thus, we deduce the desired convergence
\begin{equation*}
\lim_{n\rightarrow \infty }\limsup_{m\rightarrow \infty }\mathbb{P}\left\{
\sup_{t\in \lbrack 0,T]}\left\vert \sum_{j=1}^{[nt]}\int_{I_{j}^{n}}\left(
f(s)-f(u_{j-1})\right) dg_{m}(s)\right\vert >\varepsilon \right\} =0,
\end{equation*}%
since $\alpha +\beta >1$. This would complete the proof of our central limit
theorem.

Unfortunately, there are only few results about the asymptotic behaviour of $%
\upsilon _{\frac{1}{\beta }}(g_m)$ (the latter denotes strong $1/\beta$%
-variation on the interval $[0,1]$). Below, we shall mention some of them.
Denote by $\mathcal{W}_{p}[0,1]$ the space of functions on $[0,1]$ such that
$\upsilon _{p}<\infty $ $(p\geq 1)$ with the norm $||\cdot ||_{[p]}:=\left(
\upsilon _{p}(\cdot )\right) ^{1/p}+||\cdot ||_{\infty }.$

\begin{proposition}
(Norvai\v{s}a-Ra\v{c}kauskas, 2008, \cite{NorRac08}). Let $(X_{i})_{i=1}^{m}$
be an iid sequence, set $S_{n}(t):=\sum_{i=1}^{[mt]}X_{i},$ $t\in \lbrack
0,1].$ Then, for $p>2,$%
\begin{equation*}
\frac{1}{\sqrt{m}}S_{m}\underset{m\rightarrow \infty }{\overset{\mathcal{L}}{%
\rightarrow }}w,
\end{equation*}%
in $\mathcal{W}_{p}[0,1]$ iff $\mathbb{E}(X_{1})=0$ and $\mathbb{E}%
(X_{1}^{2})=1.$
\end{proposition}

As a consequence, if the random variables $\xi_{i,m}$ are iid with
\begin{equation*}
\mathbb{E}\left( \xi_{i,m} \right)=0, \qquad \mathbb{E}\left(
\xi_{i,m}^2\right)=\frac 1m,
\end{equation*}
we immediately deduce that $\upsilon _{p}(g_m)\underset{m\rightarrow \infty }%
{\overset{\mathcal{L}}{\rightarrow }}\upsilon _{p}(w).$

Another result that can help is the following one.

\begin{theorem}
(L\'{e}pingle 1976, \cite{Lep76}). If $p>2$ there exists a positive constant
$C$ depending on $p$, such that
\begin{equation*}
\mathbb{E}(\upsilon _{p}(M)^{1/p})\leq C \mathbb{E}(||M||_{\infty })
\end{equation*}%
for all martingales $M$.
\end{theorem}

Assume that $g_m$ is \ a martingale for a fixed $m$, which is equivalent to
say that $\{ \xi_{i,m}, 1\le i \le [mT]\}$ is a martingale difference. Then
we have%
\begin{equation*}
\mathbb{E}(\upsilon _{p}(g_m)^{1/p})\leq C \mathbb{E}(||g_m||_{\infty }),
\end{equation*}%
and if%
\begin{equation*}
g_m \overset{\mathcal{L}}{\underset{m\rightarrow \infty }{\rightarrow }}w
\end{equation*}%
in $D[0,T]$, we also have that
\begin{equation*}
\left\Vert g_m\right\Vert _{\infty }\overset{\mathcal{L}}{\underset{%
m\rightarrow \infty }{\rightarrow }}\left\Vert w\right\Vert _{\infty }.
\end{equation*}
Using the tightness condition (\ref{tight}) on $g_m$ and Doob's inequality,
we obtain that
\begin{equation*}
\mathbb{E}(||g_m||_{\infty })<C.
\end{equation*}%
Now, by the Skorohod representation theorem and dominated convergence
theorem, we have that
\begin{equation*}
\limsup_{m\rightarrow \infty} \mathbb{E}(\upsilon_{p}^j(g_m)^{1/p})\leq
Cn^{-\beta }||w||_{\beta }.
\end{equation*}%
and we can obtain the central limit theorem.

So, in the both cases mentioned above we need additional conditions on the
process $g_m$, in order to get the desired result. However, other
interesting examples are not covered by above methods. For instance,
consider a fractional Brownian motion $(B^H_t)_{t\in [0,T]}$ with Hurst
parameter $H\in (0,3/4)$ and define
\begin{equation*}
g_m (t)= \frac{1}{\sqrt{m}} \sum_{i=1}^{[mt]} \left(m^{2H}(B_{t_i}^H -
B_{t_{i-1}}^H)^2 -1 \right).
\end{equation*}
It is well known that $g_m \overset{\mathcal{L}}{\rightarrow} w$ on $D[0,T]$%
, but, to the best of our knowledge, there is less known about the
asymptotic behaviour of the strong $p$-variation of $g_m$. For this reason
we will develop a new technique, which does not rely on $p$-variation
concepts.

\section{Proof of the convergence (\protect\ref{vancov})}

In this section we are going to prove (\ref{vancov}) and, therefore,
complete the proof of Theorem \ref{thm1}, using techniques of fractional
calculus. We refer the reader to \cite{SamKilMar93} for a detailed
exposition of fractional calculus.

\bigskip
\begin{proof}
Fix $\gamma \in (0,1)$ such that $1/2<\gamma <\alpha $. Throughout the proof
all positive constants are denoted by $C$, although they may change from
line to line. Denote by $\beta _{j}$ the smallest integer greater or equal
than $mu_{j}$. Then, an integer $i$ belongs to $I_{n}(j)$ if and only if $%
u_{j-1}\leq \frac{i}{m}<\frac{\beta _{j}}{m}$. Let $J_{j}^{n,m}$ be the
interval $\left[ u_{j-1},\frac{\beta _{j}}{m}\right) $. With this notation
we can write
\begin{equation*}
R_{n,m}(t)=\sum_{j=1}^{[nt]}\sum_{i\in I_{n}(j)}\left(
f(t_{i})-f(u_{j-1})\right) \xi
_{i,m}=\sum_{j=1}^{[nt]}\int_{J_{j}^{n,m}}(f(s)-f(u_{j-1}))\mathrm{d}%
g_{m}(s).
\end{equation*}%
Set
\begin{equation*}
R_{n,m,j}:=\int_{J_{j}^{n,m}}(f(s)-f(u_{j-1}))\mathrm{d}g_{m}(s).
\end{equation*}%
We have that for any $0\leq a<b\leq T$, the identity
\begin{equation}
\int_{\lbrack a,b)}(f(s)-f(a))\mathrm{d}g_{m}(s)=\int_{a}^{b}D_{a+}^{\gamma
}f_{a}(s)D_{b-}^{1-\gamma }(g_{m})_{b-}(s)\mathrm{d}s  \label{formula}
\end{equation}%
holds, where $f_{a}(s)=f(s)-f(a),(g_{m})_{b-}(s)=g_{m}(s)-g_{m}(b-)$ for $%
s\in \lbrack a,b)$ and $D_{a+}^{\gamma }$ and $D_{b-}^{1-\gamma }$ are the
fractional derivative operators defined by
\begin{equation*}
D_{a+}^{\gamma }f_{a}(s)=\frac{1}{\Gamma \left( 1-\gamma \right) }\left(
\frac{f(s)-f(a)}{(s-a)^{\gamma }}+\gamma \int_{a}^{s}\frac{f(s)-f(y)}{%
(s-y)^{\gamma +1}}\mathrm{d}y\right) ,
\end{equation*}%
and%
\begin{eqnarray*}
&&D_{b-}^{1-\gamma }(g_{m})_{b-}(s) \\
&=&\frac{1}{\Gamma \left( \gamma \right) }\left( \frac{g_{m}(s)-g_{m}(b-)}{%
(b-s)^{1-\gamma }}+(1-\gamma )\int_{s}^{b}\frac{g_{m}(s)-g_{m}(y)}{%
(y-s)^{2-\gamma }}\mathrm{d}y\right) .
\end{eqnarray*}%
Notice that these operators are well defined by the $\alpha $-H\"{o}lder
continuity of $f$, \ with $\alpha >\gamma $ and since $g_{m}$ is piecewise
constant. The identity (\ref{formula}) can be found e.g. in \cite[Theorem
3.1 (v)]{Zah}. Now, if we take $a=u_{j-1},$ we can estimate $D_{a+}^{\gamma
}f_{a}\left( s\right) $ by the H\"{o}lder norm of $f$ and in this way we
obtain that
\begin{equation*}
|D_{a+}^{\gamma }f_{a}\left( s\right) |\leq C\Vert f\Vert _{\alpha
}n^{-\alpha +\gamma }=Gn^{-\alpha +\gamma },
\end{equation*}%
for some random variable $G$. As a consequence, if we put $a=u_{j-1}$ and $b=%
\frac{\beta _{j}}{m}$, we deduce the inequality
\begin{eqnarray*}
R_{n,m,j} &:&=\left\vert \int_{a}^{b}D_{a+}^{\gamma
}f_{a}(s)D_{b-}^{1-\gamma }(g_{m})_{b-}(s)\mathrm{d}s\right\vert \\
&\leq &\frac{Gn^{-\alpha +\gamma }}{\Gamma \left( 1-\gamma \right) }%
\int_{u_{j-1}}^{\beta _{j}/m}|D_{b-}^{1-\gamma }(g_{m})_{b-}(s)|\mathrm{d}s
\\
&\leq &\frac{Gn^{-\alpha +\gamma }}{\Gamma \left( 1-\gamma \right) }%
\sum_{k=\beta _{j-1}}^{\beta _{j}}\int_{t_{k-1}}^{t_{k}}|D_{b-}^{1-\gamma
}(g_{m})_{b-}(s)|\mathrm{d}s,
\end{eqnarray*}%
The last inequality follows from the inclusion $\left[ u_{j-1},\frac{\beta
_{j}}{m}\right] \subset \left[ t_{\beta _{j-1}-1},t_{\beta _{j}}\right] $.
Notice also that $g_{m}(b-)=g_{m}(t_{\beta _{j}-1})$.

Suppose that $t_{k-1}<s<t_{k}$. Then we obtain the identity,
\begin{eqnarray*}
&&D_{b-}^{1-\gamma }(g_{m})_{b-}(s) \\
&=&\frac{1}{\Gamma \left( \gamma \right) }\left( \frac{%
g_{m}(t_{k-1})-g_{m}(t_{\beta _{j}-1})}{\left( t_{\beta _{j}}-s\right)
^{1-\gamma }}+(1-\gamma )\int_{s}^{t_{\beta _{j}}}\frac{g_{m}(t_{k-1})-g_{m}%
\left( y\right) }{\left( y-s\right) ^{2-\gamma }}dy\right)  \\
&=&\frac{1}{\Gamma \left( \gamma \right) }\left( \frac{%
g_{m}(t_{k-1})-g_{m}(t_{\beta _{j}-1})}{\left( t_{\beta _{j}}-s\right)
^{1-\gamma }}+(1-\gamma )\sum_{l=k+1}^{\beta _{j}}\int_{t_{l-1}}^{t_{l}}%
\frac{g_{m}(t_{k-1})-g_{m}\left( t_{l-1}\right) }{\left( y-s\right)
^{2-\gamma }}dy\right)  \\
&=&\frac{1}{\Gamma \left( \gamma \right) }\left( \frac{%
g_{m}(t_{k-1})-g_{m}(t_{\beta _{j}-1})}{\left( t_{\beta _{j}}-s\right)
^{1-\gamma }}\right.  \\
&&\left. -\sum_{l=k+1}^{\beta
_{j}}(g_{m}(t_{k-1})-g_{m}(t_{l-1}))[(t_{l}-s)^{\gamma
-1}-(t_{l-1}-s)^{\gamma -1}]\right) .
\end{eqnarray*}%
Therefore,
\begin{eqnarray*}
&&|D_{b-}^{1-\gamma }(g_{m})_{b-}(s)| \\
&\leq &\frac{1}{\Gamma (\gamma )}\left( \sum_{l=k+1}^{\beta
_{j}}|g_{m}(t_{k-1})-g_{m}(t_{l-1})|[(t_{l-1}-s)^{\gamma
-1}-(t_{l}-s)^{\gamma -1}]\right.  \\
&&\left. +|g_{m}(t_{k}-1)-g_{m}(t_{\beta _{j}-1})|\left( t_{\beta
_{j}}-s\right) ^{\gamma -1}\right) .
\end{eqnarray*}%
Integrating in the variable $s$ yields,
\begin{eqnarray*}
&&\int_{t_{k-1}}^{t_{k}}|D_{b-}^{1-\gamma }(g_{m})_{b-}(s)|ds \\
&\leq &Cm^{-\gamma }\sum_{l=k+1}^{\beta
_{j}}|g_{m}(t_{k-1})-g_{m}(t_{l-1})||(l-1-k)^{\gamma }-2(l-k)^{\gamma
}+(l-k+1)^{\gamma }| \\
&&+Cm^{-\gamma }|g_{m}(t_{k-1})-g_{m}(t_{\beta _{j}-1})|[(\beta
_{j}-k+1)^{\gamma }-(\beta _{j}-k)^{\gamma }].
\end{eqnarray*}%
Then, for any constant $K>0$, and $\varepsilon >0$
\begin{equation*}
\mathbb{P}(\sup_{t\in \lbrack 0,T]}\left\vert R_{n,m}(t)\right\vert
>\varepsilon )\leq \mathbb{P(}G>K)+\frac{1}{\varepsilon }\mathbb{E}\left(
\sup_{t\in \lbrack 0,T]}\left\vert R_{n,m}(t)\right\vert \mathbf{1}_{\{G\leq
K\}}\right) .
\end{equation*}%
Now, we have
\begin{eqnarray*}
&&\mathbb{E}\left( \sup_{t\in \lbrack 0,T]}\left\vert R_{n,m}(t)\right\vert
\mathbf{1}_{\{G\leq K\}}\right) \leq \sum_{j=1}^{[nT]}\mathbb{E}\left(
|R_{n,m,j}|\mathbf{1}_{\{G\leq K\}}\right)  \\
&&\qquad \leq KCn^{-\alpha +\gamma }m^{-\gamma
}\sum_{j=1}^{[nT]}\sum_{k=\beta _{j-1}}^{\beta _{j}}\sum_{l=k+1}^{\beta _{j}}%
\mathbb{E}\left( |g_{m}(t_{k-1})-g_{m}(t_{l-1})|\right)  \\
&&\qquad \qquad \times |(l-1-k)^{\gamma }-2(l-k)^{\gamma }+(l-k+1)^{\gamma }|
\\
&&\qquad \qquad +KCn^{-\alpha +\gamma }m^{-\gamma
}\sum_{j=1}^{[nT]}\sum_{k=\beta _{j-1}}^{\beta _{j}}\mathbb{E}\left(
|g_{m}(t_{k-1})-g_{m}(t_{\beta _{j}-1})|\right) [(\beta _{j}-k+1)^{\gamma
}-(\beta _{j}-k)^{\gamma }].
\end{eqnarray*}%
Due to tightness condition \eqref{tight} we obtain that
\begin{equation*}
\mathbb{E}(|g_{m}(t_{k-1})-g_{m}(t_{l-1})|)\leq C|t_{k}-t_{l}|^{1/2},
\end{equation*}%
hence,
\begin{eqnarray*}
&&\mathbb{E}\left( \sup_{t\in \lbrack 0,T]}\left\vert R_{n,m}(t)\right\vert
\mathbf{1}_{\{G\leq K\}}\right) \leq KCn^{-\alpha +\gamma }m^{-\gamma
}\sum_{j=1}^{[nT]}\sum_{k=\beta _{j-1}}^{\beta _{j}}\sum_{l=k+1}^{\beta _{j}}%
\sqrt{l-k} \\
&&\qquad \times |(l-1-k)^{\gamma }-2(l-k)^{\gamma }+(l-k+1)^{\gamma }| \\
&&\qquad +KCn^{-\alpha +\gamma }m^{-\gamma }\sum_{j=1}^{[nT]}\sum_{k=\beta
_{j-1}}^{\beta _{j}}\sqrt{\beta _{j}-k}[(\beta _{j}-k+1)^{\gamma }-(\beta
_{j}-k)^{\gamma }] \\
&&\qquad \qquad \leq KCn^{-\alpha +\gamma }m^{-\gamma -\frac{1}{2}%
}\sum_{j=1}^{[nT]}\left( \left( \beta _{j}-\beta _{j-1}\right)
\sum_{k=1}^{\beta _{j}-\beta _{j-1}}k^{\gamma -\frac{3}{2}%
}+\sum_{k=1}^{\beta _{j}-\beta _{j-1}}k^{\gamma -\frac{1}{2}}\right)  \\
&&\qquad \qquad \leq KCn^{-\alpha +\gamma }m^{-\gamma -\frac{1}{2}%
}\sum_{j=1}^{[nT]}\left( \beta _{j}-\beta _{j-1}\right) ^{\gamma +\frac{1}{2}%
}.
\end{eqnarray*}%
Then, since $\beta _{j}-\beta _{j-1}\leq \frac{m}{n}+1$, we conclude that
\begin{equation*}
\mathbb{E}(\sup_{t\in \lbrack 0,T]}\left\vert R_{n,m}(t)\right\vert \mathbf{1%
}_{\{G\leq K\}})\leq KCn^{1/2-\alpha }\rightarrow 0,
\end{equation*}%
since $\alpha >1/2.$ Therefore, by letting $K$ goes to infinity, \ we obtain
\begin{equation*}
\lim_{n\rightarrow \infty }\lim_{m\rightarrow \infty }\mathbb{P(}\sup_{t\in
\lbrack 0,T]}\left\vert R_{n,m}(t)\right\vert >\varepsilon )=0.
\end{equation*}
\end{proof}

\end{document}